\documentclass[12pt]{article}               % The use of LaTeX2e is preferred.
\usepackage{geometry}
\usepackage{cite}
\usepackage{multirow}
\usepackage{amsmath}
\usepackage{amssymb}
\usepackage[usenames]{color}
\usepackage{graphicx}          % Include this line if your
 \usepackage{amsthm}                              % document contains figures,
\usepackage{epsfig}    % or this line, depending on which
 \usepackage{helvet}

\newtheorem{theorem}{Theorem}

\newtheorem{proposition}{Proposition}

\title{\large On the Eigenvalue Distribution for a Beam with Attached Masses}

\author{Julia Kalosha\thanks{
Institute of Applied Mathematics and Mechanics, National Academy of Sciences of Ukraine, Slovyansk ({\tt\small  julykucher@gmail.com})
        \newline
$^{**}$Max Planck Institute for Dynamics of Complex Technical Systems, Magdeburg, Germany
({\tt\small zuyev@mpi-magdeburg.mpg.de, benner@mpi-magdeburg.mpg.de})\newline
} , Alexander Zuyev$^{*,**}$, and Peter Benner$^{**}$}
\date{}

\begin{document}

\maketitle
\thispagestyle{empty}
%\pagestyle{empty}

%%%%%%%%%%%%%%%%%%%%%%%%%%%%%%%%%%%%%%%%%%%%%%%%%%%%%%%%%%%%%%%%%%%%%%%%%%%%%%%%
\begin{abstract}
We study a mathematical model of a hinged flexible beam with piezoelectric actuators and electromagnetic shaker in this paper.
The shaker is modelled as a mass and spring system attached to the beam. To analyze free vibrations of this mechanical system, we consider the corresponding spectral problem for a fourth-order differential operator with interface conditions that characterize the shaker dynamics.
The characteristic equation is studied analytically, and asymptotic estimates of eigenvalues are obtained.
The eigenvalue distribution is also illustrated by numerical simulations under a realistic choice of mechanical parameters.
\end{abstract}

\section{Introduction}\label{section_intro}
The stabilization problem for distributed parameter mechanical systems with flexible beams has been the subject of investigations of many authors. Without any pretence to represent an exhaustive list, we just mention some of the most related works.
In the monograph~\cite{LGM1999}, the problem of strong asymptotic stabilization of vibration systems with Euler--Bernoulli beams has been investigated.
In~\cite{O2000}, spectral methods for the investigation of strong stability of distributed parameter systems in Hilbert spaces are presented on a common basis, and some examples of controlled beams are analyzed.
The stabilization of an essentially nonlinear system consisting of a beam attached to the disc with torque control is achieved in~\cite{C2007}.

The dynamic behavior of a non-uniform cantilever Euler--Bernoulli beam with boundary control is investigated in the paper~\cite{G2002}.
It is emphasized that asymptotic properties of the eigenvalues of the corresponding infinitesimal generator play an important role for proving the exponential stability of the closed-loop system. A numerical scheme that preserves the exponential stability of the original continuous Euler--Bernoulli model with a clamped end and boundary control is proposed in~\cite{LG2019} without introducing additional numerical viscosity.
The Euler--Bernoulli beam with one clamped end is considered in~\cite{S2019} as a control system with two-dimensional input and output.
It is shown that the spectrum of the generator is located in the open upper half-plane for the considered model, and asymptotic approximations of the eigenvalues are proposed. In the paper~\cite{KS2019}, the general eigenvalue problem is formulated for the Timoshenko and the Euler--Bernoulli beam models under different boundary conditions.
A regularization process for the characteristic equation is discussed, and numerically stable forms of the original hyperbolic expressions are obtained.

Over the last few years, stability and stabilization problems for laminated beams and beam networks have received a lot of attention due to their rich theoretical context and engineering applications (see, e.g.,~\cite{LLLL2019,ARN2019,A2019,M2018} and references therein).
Stability conditions for elastic structures formed by serially connected flexible beams with different types of boundary conditions are obtained and illustrated by numerical examples in~\cite{CDKP1987,LLS1994}.
In the paper~\cite{SPI2019}, the covariance functions are computed for the transverse displacement of the Euler--Bernoulli beam under stochastic excitation.
The Euler--Bernoulli model with one clamped end and another simply supported end under the action of additive white noise is studied in~\cite{CGCJF2019}. It is shown that there exists a global random attractor for the considered stochastic system, and the Hausdorff dimension of the global random attractors is estimated.

In contrast to the above-mentioned publications,  in this work we focus on the mathematical model of a simply supported beam of length $l$ with $k$ piezoelectric actuators and an electromagnetic shaker described in~\cite{ZK2013}.
This mathematical model is represented by the following abstract differential equation~\cite{ZK2013}:
\begin{equation}
\frac d{dt}\;\xi(t)={\cal A} \xi(t)+{\cal B} y,\quad \xi(t)\in {X},\;y\in\mathbb R^{k+1}.
\label{infdim}
\end{equation}
Here
$$
X = \{\xi=(u,v,p,q)^T\,|\,u\in \overset{\;\,\circ}H\,\!^2(0,l),\; v\in L^2(0,l),\; p,q\in\mathbb C\}
$$
is the Hilbert space with the inner product
\begin{equation*}
\begin{aligned}
\left\langle\left(
                                     \begin{array}{c}
                                       u_1 \\
                                       v_1 \\
                                       p_1 \\
                                       q_1 \\
                                     \end{array}
                                   \right),\left(
                                             \begin{array}{c}
                                               u_2 \\
                                               v_2 \\
                                               p_2 \\
                                               q_2 \\
                                             \end{array}
                                           \right)
\right\rangle_{X}=\int\limits_0^l\left(E(x)I(x)u_1''(x) \bar u_2''(x)\right.&+\left.\rho(x)v_1(x) \bar v_2(x)\right)dx\\
&+\varkappa p_1 \bar p_2+mq_1\bar q_2,
\end{aligned}
\end{equation*}
and $y=(F,M_1,..., M_k)^T$ is considered as the control. The components of $y$ have the following physical meaning~\cite{ZK2013}: $F$ is the control force applied by the shaker to the beam at a point $x=l_0$, $l_0\in (0,l)$, and $M_j$ characterizes the action of the $j$-th piezoelectric actuator, $j=1,2,...,k$.
Here $\rho(x)>0$, ${E(x)I(x)\in C^2[0,l]}$, $E(x)I(x)>0$ for all $x\in(0,l)$, $m$ and $\varkappa$ are assumed to be positive constants.
Throughout the text, the prime stands for the derivative with respect to the spatial variable $x$.
The Sobolev spaces $\overset{\;\,\circ}H\,\!^k(0,l)$ consist of all functions ${u\in H^k(0,l)}$ such that ${u(0)=u(l)=0}$.

The linear differential operator $\cal A$ with the domain
\begin{equation*}
{\cal D}({\cal A})=\left\{\xi=\left(\begin{array}{c}
                 u \\
                 v \\
                 p \\
                 q \\
               \end{array}\right)\in {X}:\quad
\begin{aligned}&u\in H^4(0,l_0)\cap H^4(l_0,l),\\
&u''(0)=u''(l)=0,\\
&u''(l_0-0)=u''(l_0+0),\\
&v\in\overset{\;\,\circ}H\,\!^2(0,l),\\
&p=u(l_0),\quad q=v(l_0)
\end{aligned}
\right\}
\end{equation*}
is defined by the following rule:
\begin{equation}
{\cal A}:\xi=\left(
               \begin{array}{c}
                 u \\
                 v \\
                 p \\
                 q \\
               \end{array}
             \right)
             \mapsto {\cal A}\xi=
             \left(
               \begin{array}{c}
                 v \\
                 -\frac1\rho(EIu'')'' \\
                 q \\
                 -\frac1m\left\{\varkappa p-(EIu'')'\Big|_{l_0-0}+(EIu'')'\Big|_{l_0+0}\right\} \\
               \end{array}
             \right).
\label{A}
\end{equation}
The linear operator ${\cal B}:\mathbb R^{k+1}\to X$ is given by its matrix
\begin{equation}
{\cal B}=\left(     \begin{array}{cccc}
               0 & 0 & \ldots & 0 \\
               0 & \frac1\rho\psi_1'' & \ldots & \frac1\rho\psi_k'' \\
               0 & 0 & \ldots & 0 \\
               \frac1m & 0 & \ldots & 0 \\
             \end{array}
           \right),
\label{B}
\end{equation}
where $\psi_j\in C^2(0,l)$ is the shape function of the $j$-th piezoelectric actuator such that ${{\rm supp}\,\psi_j\cap\{0,l_0,l\}=\emptyset}$, $j=1,2,...,k$.
A state feedback law $y={\cal K}\xi$ has been proposed in~\cite{ZK2013} such that ${\cal K}: X\to \mathbb R^{k+1}$ is a bounded linear operator,
\begin{equation}
\xi=\begin{pmatrix}u \\ v \\ p \\ q\end{pmatrix} \in X\mapsto {\cal K}\xi =-\begin{pmatrix}\alpha_0 q
\\ \alpha_1 \int_0^l \psi_1''(x)v(x)dx \\ \vdots \\
\alpha_k \int_0^l \psi_k''(x)v(x)dx
\end{pmatrix} \in \mathbb R^{k+1},
\label{K}
\end{equation}
and $\alpha_0$, $\alpha_1$, ..., $\alpha_k$ are positive tuning parameters.
We refer to the following stability result for the closed-loop system~\eqref{infdim} with $y={\cal K}\xi(t)$.

\begin{theorem}\label{thm1}{\rm \cite{ZK2013}}
Let $\tilde {\cal A} = {\cal A}+{\cal B}{\cal K}$, $\tilde {\cal A}:{\cal D}({\cal A})\to X$, where the operators ${\cal A}$, $\cal B$, and $\cal K$ are given by~\eqref{A}--\eqref{K}, and let $\alpha_j>0$ for $j=0,1,...,k$.
Then the abstract Cauchy problem
\begin{align}\label{eq_CauchyPrb_eq}
\frac d{dt}\xi(t)&=\tilde {\cal A}\xi(t),\\\label{eq_CauchyPrb_IC}
\xi(0)&=\xi_0\in X%\\\label{eq_CauchyPrb_Control}
\end{align}
is well-posed on $t\ge 0$, and the solution $\xi = 0$ of the closed-loop system~\eqref{eq_CauchyPrb_eq} is stable in the sense of Lyapunov.
\end{theorem}

Note that the above theorem does not guarantee asymptotic stability of the equilibrium $\xi=0$ as the limit behavior of trajectories
is determined by invariant subsets of the set $M=\{\xi\in {\cal D}({\cal A})\,|\,\dot { E}(\xi)=0\}$, where $E(\xi)$ is the weak Lyapunov functional constructed in~\cite{ZK2013}.
As it is well-known, the analysis of such invariant subsets is related to spectral properties of the infinitesimal generator (see, e.g,~\cite{Z2015} and references therein).

Moreover, the distribution of eigenvalues is one of the crucial characteristics of distributed parameter systems related to their controllability, observability, and stabilizability properties because of the Ingham-type theorems~\cite{B2002} with important applications to one-dimensional wave-like equations~\cite{K1992,KS2002}.
Thus, the goal of our present paper is to study the spectral problem for $\cal A$ and analyze the distribution of roots of the corresponding frequency equation. In Section~\ref{section_spectral}, we derive an approximate frequency equation and prove its equivalence to the exact one in the sense of limit behavior of the roots. An estimate of the growth rate of eigenfrequencies is obtained in Section~\ref{section_freq}. We will present the results of numerical simulations to illustrate the behavior of the roots of both equations.
For the sake of simplicity, we will assume $E$ and $I$ to be positive constants in the sequel.

\section{Spectral Problem}\label{section_spectral}
In order to obtain the characteristic equation, we consider the spectral problem
\begin{equation}\label{eq_SpecPrb_im}
\cal A \xi=\lambda\xi,\quad\xi\in {\cal D}({\cal A}),
\end{equation}
where $\lambda=i\omega$. This problem reduces to the following set of equations with respect to the components of $\xi$:
\begin{equation}\label{eq_SpecPrb_im_cw}
\begin{aligned}
v&=\lambda u,\\
\frac{EI}\rho\frac{d^4u}{dx^4}+\lambda v&=0,\\
q&=\lambda p,\\
L-\varkappa p&=\lambda mq,
\end{aligned}
\end{equation}
and the condition $\xi\in {\cal D}({\cal A})$ yields
\begin{equation}\label{eq_BC_im}
\begin{array}{ll}
  u(0)=u(l)=0, &  \\
  u''(0)=u''(l)=0, & u^{(j)}(l_0-0)-u^{(j)}(l_0+0)=0,\quad j=\overline{0,2}, \\
  v(0)=v(l)=0, & EI(u'''(l_0-0)-u'''(l_0+0))=L, \\
  p=u(l_0),\quad q=v(l_0). &
\end{array}
\end{equation}

Let us first consider the second equation in~\eqref{eq_SpecPrb_im_cw} with respect to $u(x)$:
\begin{equation}\label{eq_specPrb_im_deq}
\frac{d^4u}{dx^4}-\mu^4u=0,\quad x\neq l_0,
\end{equation}
where
$
\mu=\left(\frac\rho{EI}\omega^2\right)^{1/4}
$
is treated as the new spectral parameter of problem~\eqref{eq_SpecPrb_im},~\eqref{eq_BC_im}.
The general solution of~\eqref{eq_specPrb_im_deq} can be represented in the intervals of continuity of $u(x)$ as
\begin{equation*}
u(x)=C_1e^{-\mu x}+C_2e^{\mu x}+C_3\sin{\mu x}+C_4\cos{\mu x},\quad x\neq l_0.
\end{equation*}
Let us denote
\begin{equation*}
U(x)=\left(
       \begin{array}{c}
         u_0(x) \\
         u_1(x) \\
         u_2(x) \\
         u_3(x) \\
       \end{array}
     \right)=\left(
               \begin{array}{c}
                u(x) \\
                 u'(x) \\
                 u''(x) \\
                 u'''(x) \\
               \end{array}
             \right),
\end{equation*}
then the differential equation that defines the eigenfunctions of~\eqref{eq_SpecPrb_im},~\eqref{eq_BC_im} can be written for $x\in(0,l_0)\cup(l_0,l)$ as
\begin{equation}\label{eq_specPrb_im_deq_vec}
\frac d{dx}U(x)=MU(x),\;\; M=\left(
            \begin{array}{cccc}
              0 & 1 & 0 & 0 \\
              0 & 0 & 1 & 0 \\
              0 & 0 & 0 & 1 \\
              \mu^4 & 0 & 0 & 0 \\
            \end{array}
          \right).
\end{equation}
The general solution of~\eqref{eq_specPrb_im_deq_vec} is represented by the matrix exponential as follows:
\begin{equation}\label{eq_GenSolSpecPrb_im}
U(x)=\left\{
\begin{array}{ll}
  e^{xM}U(0), & x\in(0,l_0), \\
  e^{(x-l)M}U(l), & x\in(l_0,l),
\end{array}
\right.
\end{equation}
where
\begin{equation*}
e^{xM}=\left(
         \begin{array}{cccc}
           z_1(x) & z_2(x) & z_3(x) & z_4(x) \\
           \mu^4z_4(x) & z_1(x) & z_2(x) & z_3(x) \\
           \mu^4z_3(x) & \mu^4z_4(x) & z_1(x) & z_2(x) \\
           \mu^4z_2(x) & \mu^4z_3(x) & \mu^4z_4(x) & z_1(x) \\
         \end{array}
       \right),
\end{equation*}
\begin{align*}
&z_1(x)=\frac12(\cosh \mu x+\cos \mu x),\quad
z_2(x)=\frac1{2\mu}(\sinh \mu x+\sin \mu x),\\
&z_3(x)=\frac1{2\mu^2}(\cosh \mu x-\cos \mu x),\quad
z_4(x)=\frac1{2\mu^3}(\sinh \mu x-\sin \mu x).
\end{align*}

The constants of integration $u_0(0)$, $u_2(0)$, $u_0(l)$, $u_2(l)$ are equal to zero due to the boundary conditions of the spectral problem, while the other four can be obtained from the interface conditions at the point $x=l_0$. These conditions may be represented as the system of linear algebraic equations
$
{\cal M}\left(
            u_1(0),
            u_3(0),
            u_1(l),
            u_3(l)
        \right)^T=0,
$
whose matrix is
\begin{equation*}
{\cal M}=\left(
           \begin{array}{cccc}
             z_2(l_0) & z_4(l_0) & -z_2(l_0-l) & -z_4(l_0-l) \\
             z_1(l_0) & z_3(l_0) & -z_1(l_0-l) & -z_3(l_0-l) \\
             \mu^4z_4(l_0) & z_2(l_0) & -\mu^4z_4(l_0-l) & -z_2(l_0-l) \\
             \mu^4z_3(l_0)-\hat\mu z_2(l_0) & z_1(l_0)-\hat\mu z_4(l_0) & -\mu^4z_3(l_0-l) & -z_1(l_0-l) \\
           \end{array}
         \right),
\end{equation*}
where
$
\hat\mu=\frac\varkappa{EI}-\frac m\rho\mu^4
$.
The last equality can be obtained from the conditions
$
\varkappa p+\lambda mq=EI(u_3(l_0-0)-u_3(l_0+0)),\;q=\lambda p
$. The value of $p$ is treated as the limit of $u_0(x)$ at $x=l_0$.
Thus, we obtain the following frequency equation:
\begin{equation}\label{eq_FreqEq}
\det{\cal M}=0,
\end{equation}
where

\begin{equation*}
\begin{aligned}
\det{\cal M}= & \frac{m}{4 \mu \rho} \left\{ (\cosh \mu (l-2 l_0) - \cosh \mu l ) \sin \mu l + (\cos \mu(l-2 l_0)-\cos \mu l ) \sinh \mu l \right\}\\
&- \frac{\sin \mu l \sinh \mu l}{\mu^2}+\frac{\varkappa}{4 EI \mu^5}\left\{
(\cosh \mu l -\cosh\mu(l-2l_0))\sin \mu l \right. \\
&\left.- (\cos \mu (l-2 l_0)+\cos \mu l) \sinh \mu l \right\}.
\end{aligned}
\end{equation*}

If $\mu$ satisfies the above equation, then ${\rm rank}\,{\cal M}\leqslant3$ and it is possible to choose a non-zero set of values $u_1(0)$, $u_3(0)$, $u_1(l)$, $u_3(l)$.
Let us take, for instance $u_3(l)=1$, then the other boundary values can be found from the following system of linear algebraic equations:
\begin{equation}
{\cal M}_3 (u_1(0), u_3(0), u_1(l))^T = (z_4(l_0-l), z_3(l_0-l), z_2(l_0-l))^T,
\label{eq_BVSys1}
\end{equation}
where
$$
{\cal M}_3 = \left(
  \begin{array}{ccc}
    z_2(l_0) & z_4(l_0) & -z_2(l_0-l) \\
    z_1(l_0) & z_3(l_0) & -z_1(l_0-l) \\
    \mu^4z_4(l_0) & z_2(l_0) & -\mu^4z_4(l_0-l) \\
  \end{array}
\right).
$$
Thus, solution~\eqref{eq_GenSolSpecPrb_im} of equation~\eqref{eq_specPrb_im_deq_vec} contains the function $u(x)$ of the following form:
\begin{equation*}
u(x)=\left\{\begin{aligned}
&z_2(x)u_1(0)+z_4(x)u_3(0),\quad x\in(0,l_0),\\
&z_2(x-l)u_1(l)+z_4(x-l),\quad x\in(l_0,l),
\end{aligned}\right.
\end{equation*}
and the set of boundary values is determined uniquely by system~\eqref{eq_BVSys1}, provided that
\begin{equation}
\det {\cal M}_3 = \frac{\sinh \mu l_0  \sin \mu l + \sin \mu l_0 \sinh \mu l }{2\mu^2} \neq 0.
\label{detM3}
\end{equation}
 Then the rest of the components of $\xi$ are determined by system~\eqref{eq_SpecPrb_im_cw}. This completes the procedure of solving the spectral problem~\eqref{eq_SpecPrb_im}.

\section{Frequency Analysis}\label{section_freq}
The determinant of ${\cal M}$ in~\eqref{eq_FreqEq} admits the following asymptotic representation for ${\mu\to+\infty}$:
\begin{equation*}
\det\,{\cal M} =\frac{m e^{\mu l}}{8 \rho \mu} \bigl(\Phi_0(\mu) +o(1)\bigr),
\end{equation*}
where \begin{equation*}
\Phi_0(\mu)= 2\sin\mu(l-l_0)\sin\mu l_0-\sin\mu l.
\end{equation*}
Thus, equation~\eqref{eq_FreqEq} is equivalent to the following one if $\mu\neq 0$:
\begin{equation}\label{eq_Phi}
\Phi_0(\mu) + \Phi_1(\mu)=0,
\end{equation}
where
$$
\small
\begin{aligned}
\Phi_1(\mu) &=e^{-\mu l} \bigl\{
2\sinh\mu l \cos\mu(l-2 l_0) -2\sinh\mu l \cos\mu l -2\cosh \mu l \sin\mu l \bigr.\\
&\bigl.+2\sin\mu l \cosh\mu(l-2 l_0)+ e^{\mu l}(\cos\mu l + \sin\mu l  -\cos\mu(l-2 l_0))-\frac{8}{m\mu}\rho\sinh\mu l \sin\mu l
\bigr\} \\
&+\frac{2 \varkappa \rho e^{-\mu l}}{EI m \mu^4}\left\{( \cosh \mu l -\cosh \mu(l-2 l_0))\sin \mu l +(\cos\mu l -\cos \mu(l-2 l_0))\sinh \mu l \right\},
\end{aligned}
$$
and $\Phi_1(\mu)\to 0$ as $\mu\to+\infty$. Note that the distribution of the roots of $\Phi_0(\mu)$ is determined by the length parameters $l$ and $l_0$ only.
We will show below that the roots of~\eqref{eq_Phi} can be approximated by the roots of
\begin{equation}\label{eq_FreqEqAp}
\Phi_0(\mu)=0,
\end{equation}
provided that $\mu>0$ is large enough.

\begin{proposition}\label{le_1}
Assume that the number ${l_0}/l$ is rational.
Let $\bar\mu_1<\bar\mu_2<...$ be the positive roots of equation~\eqref{eq_FreqEqAp}.
Then, for every $\varepsilon>0$, there is an $M>0$ such that for every $\bar\mu_j>M$ there exists a unique root ${\mu_j\in I_j=(\bar\mu_j-\varepsilon;\bar\mu_j+\varepsilon)}$ of~\eqref{eq_Phi}.
\end{proposition}

It means that every root of equation~\eqref{eq_Phi} is located in an {$\varepsilon$-neighborhood} of the corresponding root of equation~\eqref{eq_FreqEqAp}, and this neighborhood does not contain any other root of~\eqref{eq_Phi}.
Besides, it is important to ensure that there are no ``large'' roots of equation~\eqref{eq_Phi} outside the $\varepsilon$-neighborhoods mentioned above,
as stated in the next proposition.

\begin{proposition}\label{le_2}
Let the assumptions of Proposition~\ref{le_1} be satisfied, and let
\begin{equation*}
S=(M;+\infty)\setminus\bigcup_{j=1}^\infty(\mu_j-\varepsilon;\mu_j+\varepsilon).
\end{equation*}
Then $\Phi_0(\mu)+\Phi_1(\mu)\neq0$ for all $\mu\in S$.
\end{proposition}

{\em Proof of Propositions~\ref{le_1} and~\ref{le_2}.}
In case when the numbers $l_0$ and $l$ are commensurable, the function $\Phi_0(\mu)$ is periodic and its zeros $0<\bar\mu_1<\bar\mu_2<...$ are simple, $\Phi_0'(\mu_j)\neq 0$.
 Without loss of generality, we assume that a number $\varepsilon>0$ is chosen small enough so that
 $\varepsilon<(\bar\mu_{j+1}-\bar\mu_j)/2$ (the intervals $I_j$ and $I_{j+1}$ are disjoint) for all $j\ge 1$, and
\begin{equation}\label{eq_lem1}
\inf_{|\mu-\bar \mu_j|\le \varepsilon, j\ge 1}|\Phi_0'(\mu)|\geqslant K
\end{equation}
with some constant $K>0$.
By construction, all positive roots of the equation $\Phi_0(\mu)=0$ are contained in the open set $I=\cup_{j\ge 1} I_j$.
Hence, $\Phi_0(\mu)\neq 0$ at each $\mu \in S_0 = [0,+\infty)\setminus I$.
By exploiting the periodicity and continuity of $\Phi_0$, we conclude from the Weierstrass extreme value theorem that there exists a $\delta>0$ such that
\begin{equation}
|\Phi_0(\mu)|\ge \delta \quad \text{for all}\;\mu \in S_0.
\label{S0}
\end{equation}
From the mean value theorem, we know that
\begin{equation}
\Phi_0(\mu)=\Phi_0'(\chi)(\mu-\bar\mu_j),
\label{mean_value}
\end{equation}
where $\chi=\bar\mu_j+\Theta(\mu-\bar\mu_j)$, $\Theta\in(0,1)$.
Then~\eqref{eq_lem1} and~\eqref{mean_value} imply that
\begin{equation}
|\Phi_0(\mu_j\pm \varepsilon)|\ge K \varepsilon\;\text{and}\; \Phi_0(\mu_j- \varepsilon)  \Phi_0(\mu_j+\varepsilon)<0 \quad \text{for all}\;j=1,2,... \, .
\label{signs}
\end{equation}

Since $\Phi_1(\mu)\to 0$ and $\Phi_1'(\mu)\to 0$ as $\mu\to +\infty$, we take an $M>0$ such that
\begin{equation}
|\Phi_1(\mu)|< \min\{K\varepsilon,\delta \}\; \text{and}\; |\Phi_1'(\mu)|<K/2 \quad\text{for all}\; \mu> M-\varepsilon.
\label{phi1_small}
\end{equation}
From~\eqref{signs} and~\eqref{phi1_small}, it follows that the continuous function $\Phi(\mu)=\Phi_0(\mu)+\Phi_1(\mu)$
has values of opposite sign at $\bar \mu_j\pm \varepsilon$, provided that $\bar\mu_j>M$. Then there exists a $\mu_j\in I_j=(\bar \mu_j - \varepsilon; \bar \mu_j + \varepsilon)$ such that $\Phi(\mu_j)=0$ by the intermediate value theorem.

The uniqueness of the root $\mu_j$ in $I_j$ can be proved by contradiction. Let $\mu_j^*\in I_j$ be such that
\begin{equation}\label{eq_lem0}
\Phi_0(\mu_j^*)=-\Phi_1(\mu_j^*),\quad \mu_j^*\neq\mu_j.
\end{equation}
On the one hand, the integral representations
\begin{equation*}
\Phi_0(\mu_j^*)=\Phi_0(\mu_j)+\int\limits_{\mu_j}^{\mu^*}\Phi_0'(\zeta)d\zeta,\quad \Phi_1(\mu_j^*)=\Phi_1(\mu_j)+\int\limits_{\mu_j}^{\mu^*}\Phi_1'(\zeta)d\zeta
\end{equation*}
together with the condition $\Phi_0(\mu_j)=-\Phi_1(\mu_j)$ yields
\begin{equation}\label{eq_lem2}
\int\limits_{\mu_j}^{\mu^*}(\Phi_0'(\zeta)+\Phi_1'(\zeta))d\zeta=0.
\end{equation}
On the other hand, inequalities~\eqref{eq_lem1} and~\eqref{phi1_small} imply
$$
\left|\int\limits_{\mu_j}^{\mu^*}(\Phi_0'(\zeta)+\Phi_1'(\zeta))d\zeta\right| \ge |{\mu_j}-{\mu^*}|\inf_{\zeta \in I_j} |\Phi_0'(\zeta)+\Phi_1'(\zeta)|\ge
$$
$$
\ge |{\mu_j}-{\mu^*}|\bigl(\inf_{\zeta \in I_j} |\Phi_0'(\zeta)| - \sup_{\zeta \in I_j} |\Phi_1'(\zeta)| \bigr) \ge \frac{K |{\mu_j}-{\mu^*}|}{2}>0.
$$
The above inequality contradicts~\eqref{eq_lem2}; hence, $\mu_j$ is the only zero of $\Phi(\mu)$ in $I_j$, which proves Proposition~\ref{le_1}.

By taking into account~\eqref{S0} and~\eqref{phi1_small}, we conclude that $|\Phi_0(\mu)|> |\Phi_1(\mu)|$ whenever $\mu\in S_0$ and $\mu>M$.
The latter implies the assertion of Proposition~\ref{le_2}.
$\square$

\section{Numerical Simulation Results}\label{section_num}
Let $0<\mu_1\le \mu_2 \le ... \le \mu_N$ be the first $N$ positive roots of equation~\eqref{eq_FreqEq}.
Then the corresponding eigenvalues $\lambda_j$ of~\eqref{eq_SpecPrb_im} and the modal frequencies $\nu_j$ in Hz are expressed as
$$
\lambda_j=i\sqrt{\frac{EI}{\rho}} \mu_j^2,\; \nu_j=\frac{\sqrt{EI}}{2\pi\sqrt{\rho}} \mu_j^2,\quad j=1,2,..., N.
$$
To test the approximation of $\mu_j$ by solutions of the truncated equation~\eqref{eq_FreqEqAp}, we also compute numerically the roots $\bar \mu_j$ of the transcendent equation $\Phi_0(\bar \mu)=0$ and the corresponding frequency parameters $\bar \nu_j=\frac{\sqrt{EI}}{2\pi\sqrt{\rho}} \bar \mu_j^2$ for $N=22$.
These numerical results are summarized in Table~\ref{tab_Eigenfrequencies} for the following realistic values of mechanical parameters:
$$
l=1.905\,\text{m},\; l_0=1.4\,\text{m},\;\rho_0 = 2700 \,\text{kg}/\text{m}^3,\; S = 2.25\cdot10^{-4} \text{m}^2,\; \rho=\rho_0 S,
$$
\begin{equation}
E=6.9\cdot10^{10}\,\text{Pa},\; I=1.6875\cdot10^{-10}\,\text{m}^4,\; m=0.1\,\text{kg},\; \varkappa=7\,\text{N/mm}.
\label{mechparameters}
\end{equation}
The above mechanical parameters correspond to the experimental setup % of the Institute of Mechanics and Mechatronics of the Vienna University of Technology
described in~\cite{DSK2014}.

\begin{table}{\caption{Modal frequencies under the choice of parameters~\eqref{mechparameters}}\label{tab_Eigenfrequencies}
\parbox[t]{.45\linewidth}
{\footnotesize
\begin{tabular}{|r||l|l||l|l|}
\hline
  $j$ & $\bar{\mu_j}$ & $\mu_j$ & $\bar{\nu_j}$ & $\nu_j$\\\hline\hline
  1 & 2.616 & 2.552 & 4.767 & 4.537\\\hline
  2 & 4.714 & 4.573 & 15.485 & 14.570\\\hline
  3 &   & 5.618 &   & 21.994\\\hline
  4 & 6.553 & 6.608 & 29.921 & 30.427\\\hline
  5 & 7.460 & 8.198 & 38.774 & 46.830\\\hline
  6 & 9.309 & 9.721 & 60.378& 65.850\\\hline
  7 & 11.407 & 11.494 & 90.661 & 92.061\\\hline
  8 & 13.013 & 13.142 & 117.997 & 120.342 \\\hline
  9 & 14.018 & 14.501 & 136.914 & 146.510\\\hline
  10 & 16.009 & 16.226 & 178.565 & 183.445\\\hline
  11 & 18.085 & 18.113 & 227.881 & 228.610\\\hline
\end{tabular}}\hspace{2em}
\parbox[t]{.45\linewidth}
{\footnotesize
\begin{tabular}{|r||l|l||l|l|}\hline
  $j$ & $\bar{\mu_j}$ & $\mu_j$ & $\bar{\nu_j}$ & $\nu_j$\\\hline\hline
  12 & 20.648 & 20.988 & 297.075 & 306.918 \\\hline
  13 & 22.711 & 22.846 & 359.376 & 363.683\\\hline
  14 & 24.732 & 24.734 & 426.184 & 426.273 \\\hline
  15 & 25.803 & 26.084 & 463.925 & 474.083\\\hline
  16 & 27.318 & 27.554 & 519.992 & 528.997 \\\hline
  17 & 29.411 & 29.497 & 602.735 & 606.252\\\hline
  18 & 31.318 & 31.325 & 683.413 & 683.723 \\\hline
  19 & 32.238 & 32.533 & 724.138 & 737.484\\\hline
  20 & 34.007 & 34.172 & 805.799 & 813.665\\\hline
  21 & 36.107 & 36.158 & 908.419 & 910.975 \\\hline
  22 & 37.814 & 37.863 & 996.296 & 998.922 \\\hline
\end{tabular}}}
\end{table}

Note that the truncated frequency equation~\eqref{eq_FreqEqAp} has an additional positive root $\bar\mu_0\approx 0.9949$ that does not correspond to any root of equation~\eqref{eq_FreqEq}, and there is the root $\mu_3$ of~\eqref{eq_FreqEq} while the equation~\eqref{eq_FreqEqAp} does not have any root in the corresponding interval.
These features are observed in low frequency range.
But for large values of~$\mu$, the agreement between solutions of the frequency equations~\eqref{eq_FreqEq} and~\eqref{eq_FreqEqAp} is quite good, as seen in Table~\ref{tab_Eigenfrequencies} and predicted by Propositions~\ref{le_1} and~\ref{le_2}.
%We can also compare these results with measurements obtained at the experimental setup of~\cite{DSK2014}.
%Some of the modal frequencies have been estimated by finding local maxima of the Fourier transform magnitude of the step response data provided by authors of~\cite{DSK2014}:
%$$
%\nu_1=4.80\, \text{Hz},\; \nu_2=14.42\,\text{Hz},\; \nu_4=30.38\,\text{Hz},\; \nu_6=65.76\,\text{Hz},\; \nu_8=120.44\,\text{Hz}.
%$$
%These experimental frequencies, indexed according to our numeration, show acceptable agreement with the simulated frequencies in Table~\ref{tab_Eigenfrequencies}.

\begin{figure}[ht]
\centering
\includegraphics[scale=0.4]{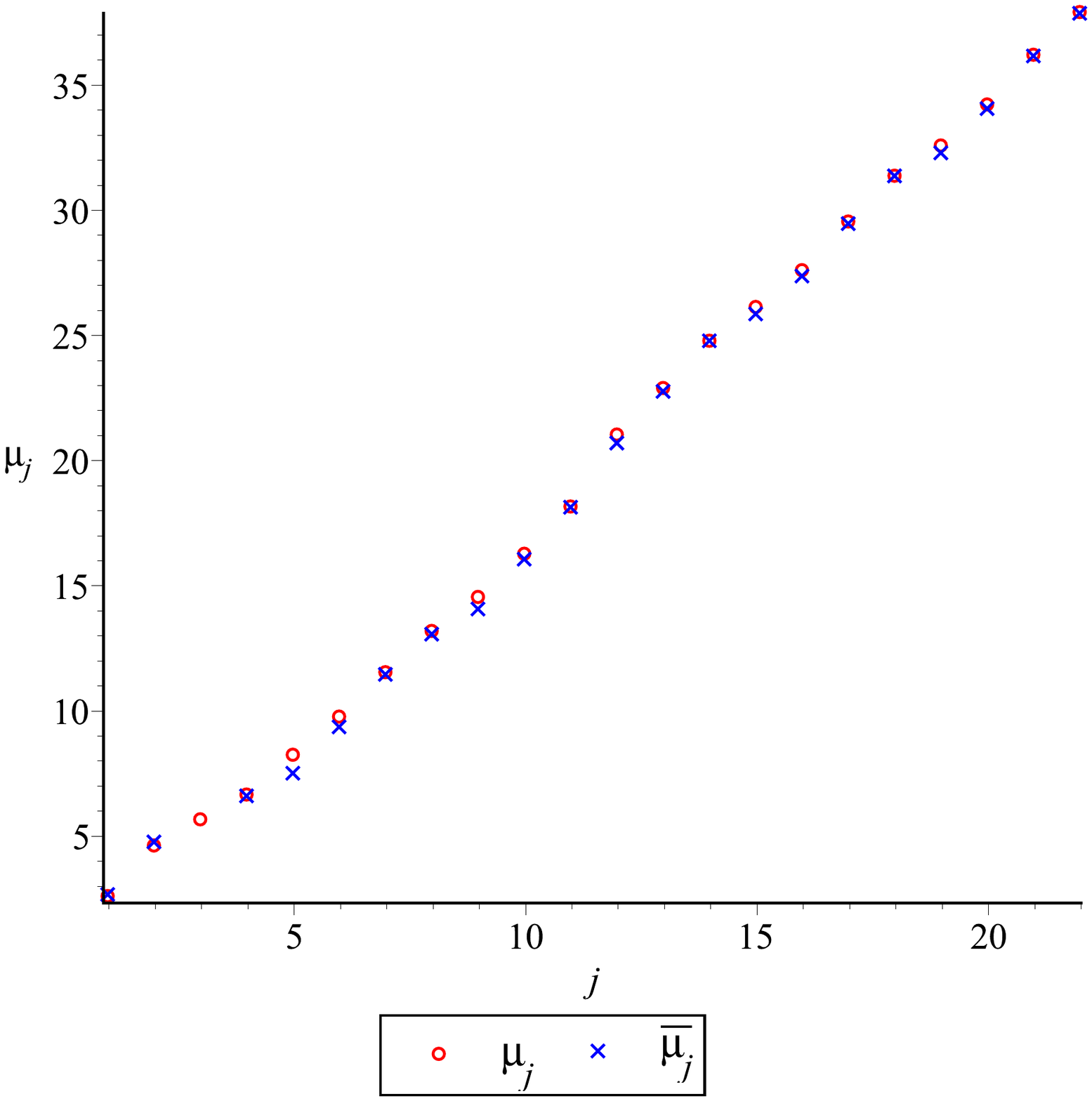}
\vskip-18mm
\caption{Spectral parameters $\mu_j$ and $\bar \mu_j$}
\label{pic_Eigenfreq}
\end{figure}

In the particular case $l=2l_0$, it is easy to write the positive roots of the transcendent equation~\eqref{eq_FreqEqAp} explicitly:
\begin{equation}
\bar \mu_{j}=\frac{\pi}l\left(\left\{\frac{j}2\right\}+2 \left[\frac{j}{2}\right]\right), \quad j=1,2,... \, .
\label{lingrowth}
\end{equation}
The above formula together with Propositions~\ref{le_1} and~\ref{le_2} implies that the spectral parameters $\mu_j$ grow linearly with $j$,
i.e. the modal frequencies $\nu_j$ and the eigenvalues $\lambda_j$ of~\eqref{eq_SpecPrb_im} grow quadratically with $j$ for large $j$ in the considered case. The growth rate of $\mu_j$ and $\bar\mu_j$ in the case of mechanical parameters~\eqref{mechparameters} is shown in Fig.~\ref{pic_Eigenfreq}.

\begin{figure}[ht]
\parbox{.45\linewidth}{\centering
\includegraphics[scale=0.275, clip]{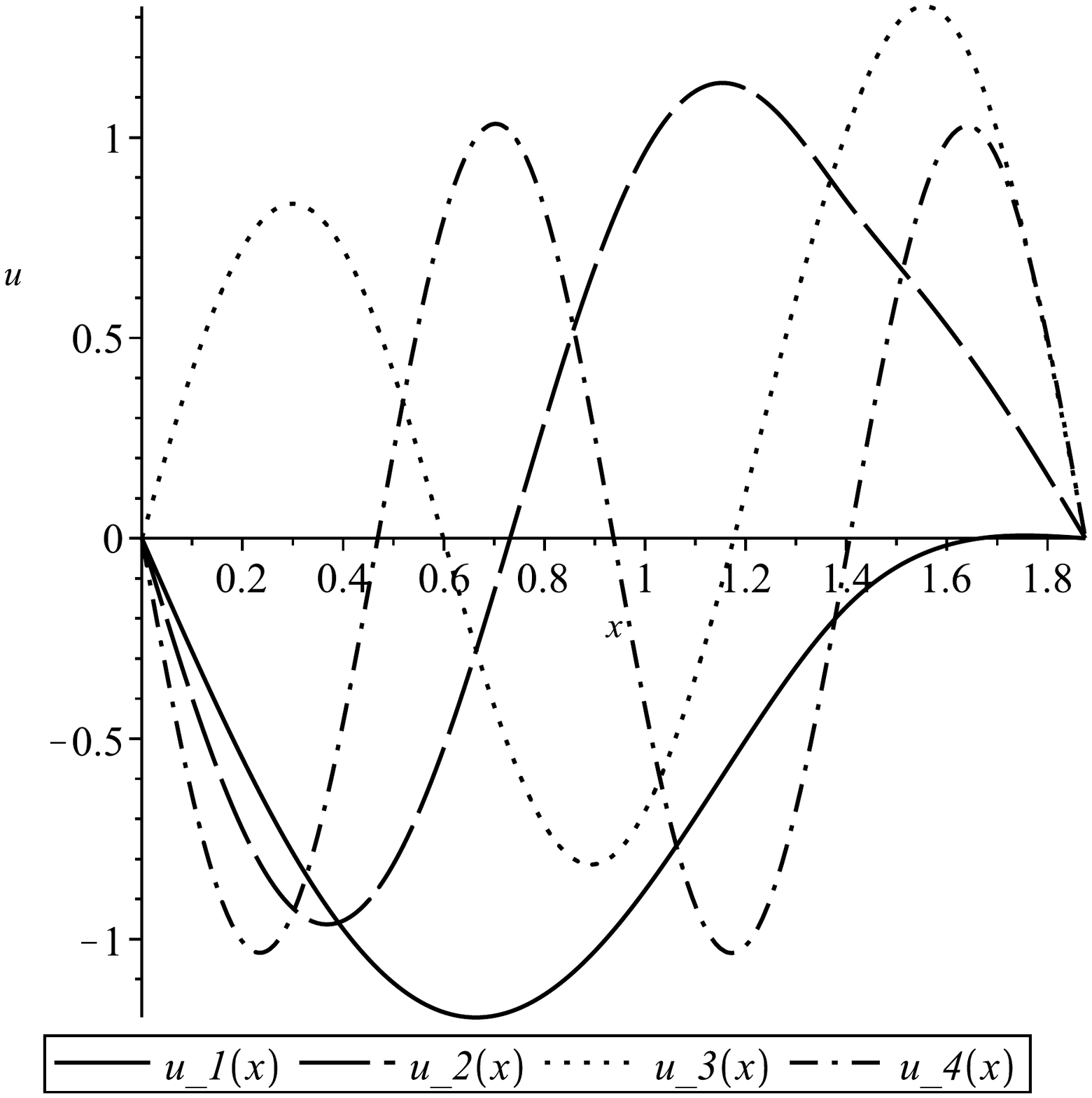}
\caption{{\footnotesize Eigenmodes for $l_0=1.4\;\text{m}$}}\label{pic_Eigenforms_lc}
}\label{pic_Eigenforms}
\parbox{.2\linewidth}{}
\parbox{.45\linewidth}{\centering
\includegraphics[scale=0.275, clip]{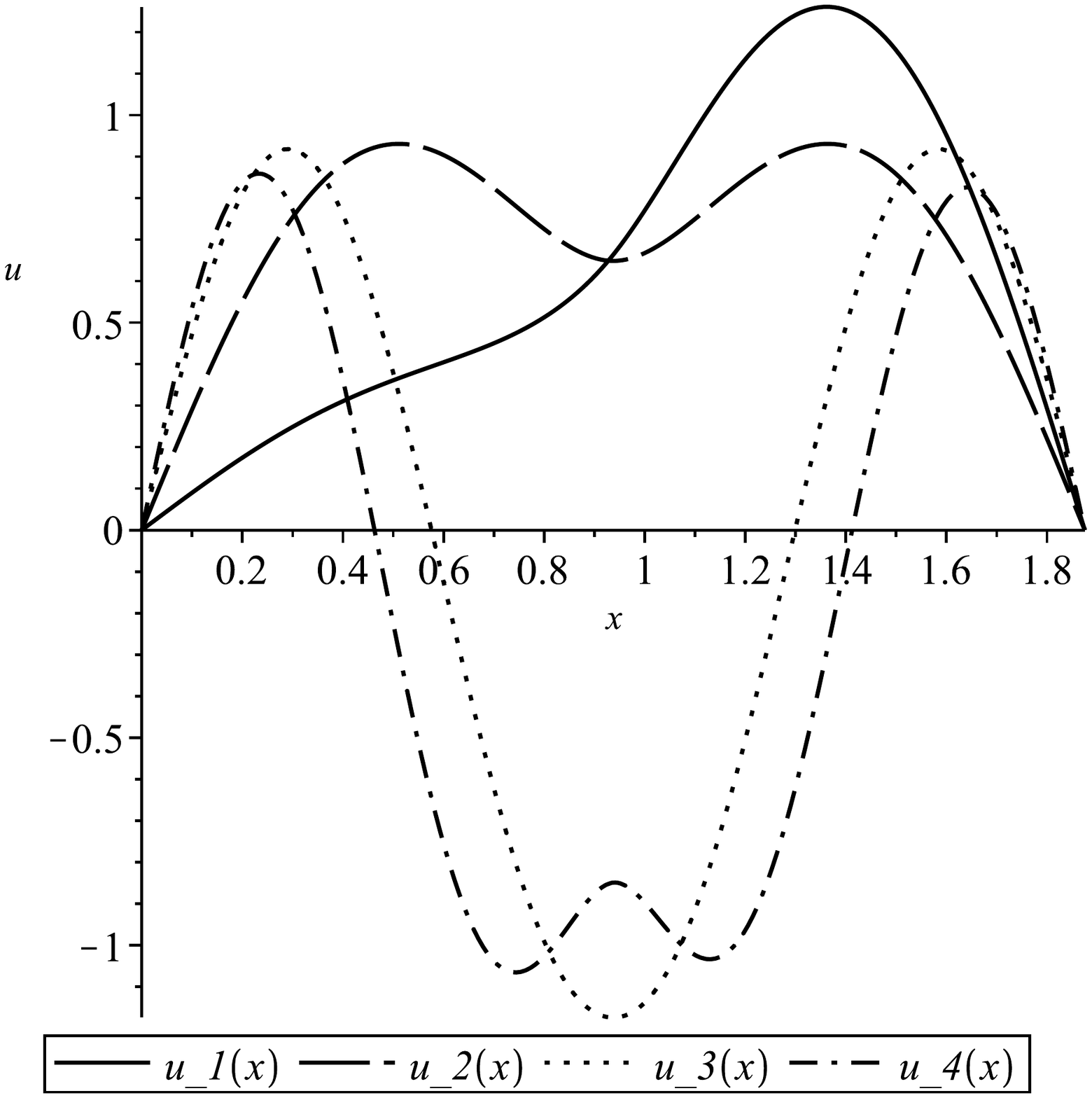}
\caption{{\footnotesize Eigenmodes for $l_0=l/2$}}}
\end{figure}

To illustrate the behavior of the corresponding eigenvectors $\xi^1$, $\xi^2$,..., $\xi^N$ of~\eqref{eq_SpecPrb_im}, we plot the graphs of their $u$-component in Figs.~\ref{pic_Eigenforms_lc} and~3.
Note that condition~\eqref{detM3} holds for each $\mu=\mu_j$ in the considered case; thus, each $\xi^j$ is uniquely defined up to normalization.
As the function $u^j(x)$ describe the transverse displacement of the beam corresponding to the spectral parameter $\mu_j$, we will refer to $u^j(x)$ as the $j$-th eigenmode of vibration.
The presented eigenmodes correspond to the choice of $l_0=1.4\;\text{m}$ in Fig.~2 and $l_0=l/2$ in Fig.~3, respectively, for $N=4$.
 In both figures, the functions $u^j(x)$ are normalized in the sense of $L^2$-norm on~$[0,l]$.
\vskip-2ex

\section{Conclusion}
The approximation results of Section~\ref{section_freq} together with the linear growth condition of the form~\eqref{lingrowth} generalize known asymptotic properties
of eigenvalues of the standard Euler--Bernoulli beam (cf.~\cite[Chapter~4]{LGM1999}) to the model with attached masses.
In future work, we expect to apply Propositions~\ref{le_1} and~\ref{le_2} for the analysis of the limit behavior of trajectories
to characterize attractors of the infinite-dimensional closed-loop system~\eqref{eq_CauchyPrb_eq} for different values of tuning parameters of the feedback law~\eqref{K}.
The results obtained in this paper are planned to be extended to other classes of elastic structures, particularly, to the rotating Timoshenko beam with attached masses~\cite{ZS2007} and beam systems with passive joints~\cite{ZS2005}.


\begin{thebibliography}{99.}%
% Use the following syntax and markup for your references if
% the subject of your book is from the field
% "Mathematics, Physics, Statistics, Computer Science"
%
% Monograph
\bibitem{LGM1999}
Luo, Z.-H., Guo, B.-Z., Morgul O.: Stability and Stabilization of Infinite Dimentional Systems with Applications. Springer, London (1999)
%
\bibitem{O2000}
Oostveen, J.: Strongly Stabilizable Distributed Parameter Systems. SIAM, Philadelphia (2000)
%
\bibitem{C2007}
Coron, J.-M.: Control and Nonlinearity. AMS, Providence (2007)
%
\bibitem{G2002}
Guo, B.-Z.: Riesz basis property and exponential stability of controlled Euler--Bernoulli beam equations with variable coefficients. SIAM Journal on Control and Optimization \textbf{40}, 1905--1923 (2002)
%
\bibitem{LG2019}
Liu, J., Guo, B.-Z.: A novel semi-discrete scheme preserving uniformly exponential stability for an Euler--Bernoulli beam. Systems \& Control Letters \textbf{134}, 1--18 (2019) doi: 10.1016/j.sysconle.2019.104518
%
\bibitem{S2019}
Shubov, M.\,A.: Location of eigenmodes of Euler-–Bernoulli beam model under fully non-dissipative boundary conditions. Proceedings of the Royal Society A \textbf{475}, 1--20 (2019) doi: 10.1098/rspa.2019.0544
%
\bibitem{KS2019}
Khasawneh, F.\,A., Segalman, D.: Exact and numerically stable expressions for Euler--Bernoulli and Timoshenko beam modes. Applied Acoustics \textbf{151}, 215--228 (2019)
%
\bibitem{LLLL2019}
Liu, W., Luan, Y., Liu, Y., Li, G.: Well‐posedness and asymptotic stability to a laminated beam in thermoelasticity of type III. Mathematical Methods in the Applied Sciences \textbf{43}, 3148--3166 (2019) doi: 10.1002/mma.6108
%
\bibitem{ARN2019}
Apalara, T.\,A., Raposo, C.\,A., Nonato, C.\,A.\,S.: Exponential stability for laminated beams with a frictional damping. Archiv der Mathematik
\textbf{114}, 471--480 (2020). doi: 10.1007/s00013-019-01427-1
%
\bibitem{A2019}
Apalara, T.\,A.: On the stability of a thermoelastic laminated beam. Acta Mathematica Scientia \textbf{39}, 1517-–1524 (2019)
%
\bibitem{M2018}
Mustafa, M.\,I.: On the stabilization of viscoelastic laminated beams with interfacial slip. Zeitschrift für angewandte Mathematik und Physik \textbf{69}. (2018) doi: 10.1007/s00033-018-0928-7.
%
\bibitem{CDKP1987}
Chen, G., Delfour, M., Krall, A., Payre, G.:. Modeling, stabilization and control of serially connected beams. SIAM Journal on Control and Optimization \textbf{25}, 526--546 (1987)
%
\bibitem{LLS1994}
Lagnese, J.E., Leugering, G., Schmidt, E.J.P.G.: Modeling, Analysis and Control of Dynamic Elastic Multi-Link Structures. Birkh{\"a}user, Boston (1994)
%
\bibitem{SPI2019}
Śniady, P., Podworna, M., Idzikowski, R.: Stochastic vibrations of the Euler--Bernoulli beam based on various versions of the gradient nonlocal elasticity theory. Probabilistic Engineering Mechanics \textbf{56}. 27--34 (2019)
%
\bibitem{CGCJF2019}
Chen, H., Guirao, J., Cao, D.Q., Jiang, J., Fan, X.: Stochastic Euler--Bernoulli beam driven by additive white noise: Global random attractors and global dynamics. Nonlinear Analysis \textbf{185}, 216--246 (2019)
%
\bibitem{ZK2013}
Zuyev, A., Kucher, J.: Stabilization of a flexible beam model with distributed and lumped controls (in Russian). Dynamical Systems. \textbf{3}(31), 25--35 (2013)
%
\bibitem{Z2015}
Zuyev, A.L.: Partial Stabilization and Control of Distributed Parameter Systems with Elastic Elements. Springer, Cham (2015)
%
\bibitem{B2002}
Baiocchi, C., Komornik, V., Loreti, P.: Ingham-Beurling type theorems with weakened gap conditions. Acta Math. Hungar. \textbf{97}, 55--95 (2002)
%
\bibitem{K1992}
Krabs, W.: On Moment Theory and Controllability of One-Dimensional Vibrating Systems and Heating Processes. Springer-Verlag, Berlin (1992)
%
\bibitem{KS2002}
Krabs, W., Sklyar, G.M.: On Controllability of Linear Vibrations. Nova Science Publishers (2002)
%
\bibitem{DSK2014}
Dullinger, C., Schirrer, A., Kozek, M.: Advanced control education: Optimal \& robust MIMO control of a flexible beam setup. IFAC Proceedings Volumes \textbf{47}(3), 9019--9025 (2014)
%
\bibitem{ZS2007}
Zuyev, A., Sawodny, O.: Stabilization and observability of a rotating Timoshenko beam model. Mathematical Problems in Engineering \textbf{2007}, 1--19 (2007) doi: 10.1155/2007/57238
%
\bibitem{ZS2005}
Zuyev, A., Sawodny, O.: Stabilization of a flexible manipulator model with passive joints. IFAC Proceedings Volumes \textbf{38}, 784--789 (2005)
\end{thebibliography}
\end{document}